\documentclass{amsart}
\usepackage{amsfonts}
\usepackage{amsmath,amssymb}
\usepackage{amsthm}
\usepackage{amscd}
\usepackage{graphics}
\usepackage{graphicx}

\usepackage[pagewise]{lineno}
\theoremstyle{definition}{
\newtheorem{Def}{{\rm Definition}}
\newtheorem{Ex}{{\rm Example}}
\newtheorem{Rem}{{\rm Remark}}
\newtheorem{Prob}{{\rm Problem}}

}
\theoremstyle{plain}
{

\newtheorem{Prop}{Proposition}
\newtheorem{Thm}{Theorem}
\newtheorem{MainThm}{Main Theorem}

}

\begin{document}
\title[On Reeb graphs of real algebraic functions which may not be planar]{Notes on Reeb graphs of real algebraic functions which may not be planar}
\author{Naoki Kitazawa}
\address{Institute of Mathematics for Industry, Kyushu University, 744 Motooka, Nishi-ku Fukuoka 819-0395, Japan\\
 TEL (Office): +81-92-802-4402 \\
 FAX (Office): +81-92-802-4405 \\
}
\email{n-kitazawa@imi.kyushu-u.ac.jp, naokikitazawa.formath@gmail.com}
\urladdr{https://naokikitazawa.github.io/NaokiKitazawa.html}
\maketitle
\begin{abstract}
The {\it Reeb graph} of a smooth function is
a graph being a natural quotient space of the manifold of the domain and the space of all connected components of preimages. Such a combinatorial and topological object roughly and compactly represents the manifold. 
Since the proposal by Sharko in 2006, reconstructing nice smooth functions and the manifolds from finite graphs in such a way that the Reeb graphs are the graphs has been important. The author has launched new studies on this, discussing construction of real algebraic
 functions. We concentrate on Reeb graphs we cannot realize as (natural) planar graphs here. Previously the graphs were planar and embedded in the plane naturally.

\end{abstract}

\section{Introduction}
\label{sec:1}

For smooth functions of some nice classes such as the class of {\rm Morse}({\rm -Bott}) {\rm functions} and a very general class in \cite{saeki}, we have their {\it Reeb graphs}. The {\it Reeb graph} of a smooth function is
a graph being the space of all connected components of preimages and as a result a quotient space of the manifold of the domain.

Such objects have been fundamental tools in understanding the manifolds roughly and compactly. \cite{reeb} is a pioneering paper.

The following problem has been first proposed by Sharko in \cite{sharko}.

\begin{Prob}
\label{prob:1}
Do we have nice smooth functions whose Reeb graphs are given graphs?
\end{Prob}

For finite graphs of suitable classes in several cases and arbitrary finite graphs, we have explicitly obtained various affirmative answers.

For example, \cite{sharko} considers some finite graphs and constructs nice smooth functions on closed surfaces. \cite{masumotosaeki} extends this to arbitrary finite graphs. \cite{michalak1} considers Morse functions on closed manifolds whose general preimages consist of spheres with suitably restricted classes of finite graphs. \cite{michalak2} is on explicit deformations of Morse functions with their Reeb graphs.

\cite{kitazawa1, kitazawa2} consider arbitrary finite graphs and first consider situations where the topologies of general preimages are as prescribed ones. \cite{saeki} is a related paper based on our informal discussions on \cite{kitazawa1}.

\begin{Prob}
Consider Problem \ref{prob:1} in the real algebraic category.
\end{Prob}

\cite{kitazawa3} is a pioneering study. See also \cite{kitazawa4, kitazawa5, kitazawa6, kitazawa7} for example.
In these studies, that the classes of graphs have been strongly restricted. 
They are in considerable cases homeomorphic to a closed interval as topological spaces. See FIGURE 1 of \cite{kitazawa3} as another case. They are very simple and planar. Furthermore, we can embed the graphs into the plane naturally.
As a related motivating study, \cite{bodinpopescupampusorea} studies graphs embedded into the plane naturally and regarded as graphs which regions in the plane surrounded by so-called {\it non-singular} real algebraic connected curves naturally collapse to and which are so-called generic graphs. We have concentrated
 on graphs of very explicit classes being subclasses of such a class of graphs.

In this note, we discuss explicit (Reeb) graphs we cannot embed into the plane in natural ways. One of our main results is as follows. Some notions and terminologies and the notation will be presented later more rigorously.

\begin{MainThm}
\label{mthm:1}
Let $D \subset {\mathbb{R}}^k$ be an open and connected set satisfying the following conditions. This respects some of Theorem \ref{thm:2}, presented later.
\begin{itemize}
\item[C1] The closure $\overline{D}$ is compact and connected. $\overline{D}-D$ is the disjoint union of finitely many {\rm non-singular} real algebraic hypersurfaces $S_j$ of dimension $k-1$ in ${\mathbb{R}}^k$, indexed by $j \in J$ in a finite set $J$ of size $l>0$, and has no boundary.
\item[C2] $D$ is, as in Definition \ref{def:4}, an {\rm NC domain}. In other words, $S_j$ is represented as some disjoint union of connected components
 of the zero set of a real polynomial $f_j$ and for some small open neighborhood $U_D \supset \overline{D}$, $D$ is represented as the intersection of {\rm (}the closure $\overline{U_D}$ of{\rm )} $U_D$ and ${\bigcap}_{j \in J} \{x \in {\mathbb{R}}^k \mid f_j(x)>0\}$
 and the closure $\overline{D}$ is represented as the intersection of {\rm (}the closure $\overline{U_D}$ of{\rm )} $U_D$ and ${\bigcap}_{j \in J} \{x \in {\mathbb{R}}^k \mid f_j(x) \geq 0\}$.
\item[C3] There exists a finite and connected graph $K_D$ enjoying the following properties.
\begin{itemize}
\item[C3.1] Its underlying space consists of all connected components of preimages considered for the restriction of the projection of ${\mathbb{R}}^k$ to the first component to the subset $\overline{D}$. It is also regarded as the quotient space of $\overline{D}$.
We also have a natural map $g_{K_D}:K_D \rightarrow \mathbb{R}$ by considering the value of the restriction of the projection of ${\mathbb{R}}^k$ to the first component at each point of $p \in K_D$, representing some preimage.
\item[C3.2] A point in the underlying space is a vertex if and only if it is a connected component containing some singular points of the natural smooth function defined as the restriction of the projection of ${\mathbb{R}}^k$ to the first component to the subset $\overline{D}-D$.
\item[C3.3] $g_{K_D}$ is a piecewise smooth function and injective on each edge of the graph $K_D$.
\end{itemize}
\end{itemize}
Suppose that two distinct values $t_1<t_2$ of the function of $g_{K_D}$ satisfy the following conditions.
\begin{enumerate}
\item[C4] ${g_{K_D}}^{-1}(t_1)$ and ${g_{K_D}}^{-1}(t_2)$ contain no vertices of $K_D$.
\item[C5] Either of the following holds.
\begin{enumerate}
\item[C5.1] ${g_{K_D}}^{-1}(t_1)$ contains two distinct points $p_{t_1,1}$ and $p_{t_1,2}$ and ${g_{K_D}}^{-1}(t_2)$ contains one point $p_{t_2}$. We have some embedded arc $e_{t_1,t_2,j}$ in the graph $K_D$ connecting $p_{t_1,j}$ and $p_{t_2}$ making the image $g_{K,D}({\rm Int}\ e_{t_1t_2,j})$ of the interior of each arc and the open interval $(t_1,t_2)$ agree for each $j$. Furthermore, the intersection $e_{t_1,t_2,1} \bigcap e_{t_1,t_2,2}$ contains a small connected embedded arc containing $p_{t_2}$.
\item[C5.2] ${g_{K_D}}^{-1}(t_1)$ contains one point $p_{t_1}$ and ${g_{K_D}}^{-1}(t_2)$ contains two distinct points $p_{t_2,1}$ and $p_{t_2,2}$. We have some embedded arc $e_{t_1,t_2,j}$ in the graph $K_D$ connecting $p_{t_1}$ and $p_{t_2,j}$ making the image $g_{K,D}({\rm Int}\ e_{t_1,t_2,j})$ of the interior of each arc and the open interval $(t_1,t_2)$ agree for each $j$. Furthermore, the intersection $e_{t_1,t_2,1} \bigcap e_{t_1,t_2,2}$ contains a small connected embedded arc containing $p_{t_1}$.
\end{enumerate}
\end{enumerate}
Then we have a family $\{K_{D,t_1,t_2,i}\}$ of graphs indexed by positive integers $i$ enjoying the following properties.
\begin{enumerate}
\item $K_{D,t_1,t_2,i_1}$ and $K_{D,t_1,t_2,i_2}$ are not isomorphic for $i_1 \neq i_2$.
\item For a sufficiently large integer $m_i>0$, we have a suitable $m_i$-dimensional non-singular real algebraic closed and connected manifold $M_i$ and a smooth real algebraic function $f_i:M_i \rightarrow \mathbb{R}$ whose Reeb graph $W_{f_i}$ is isomorphic to $K_{D,t_1,t_2,i}$.
\item For the {\rm (}uniquely defined{\rm )} function $\bar{f_i}:W_{f_i} \rightarrow \mathbb{R}$, satisfying $f_i=\bar{f_i} \circ q_{f_i}$ by the definition, it cannot be represented as the composition of any embedding into ${\mathbb{R}}^2$ with the projection to the first component.    
\end{enumerate}
\end{MainThm}

Functions on graphs such as $g_{K_D}$ here and $g_K$ in Theorem \ref{thm:2} are important in our related studies. See also Definition 2.1 of \cite{matsuzaki} for example.

The next section is for preliminary. We rigorously introduce fundamental terminologies, notions and notation we need.
The third section is on our Main Theorems. We also present additional results as Main Theorems in addition. 
All of them are on realizing graphs we cannot embed into the plane in a canonical way as the Reeb graphs of explicit real algebraic functions. Remark \ref{rem:1} with Theorem \ref{thm:2} is a remark on realizing graphs in this way for graphs we can embed into the plane in a canonical way.
After Remark \ref{rem:1}, we present additional related remarks and examples. \\
\ \\
{\bf Acknowledgement, grants and data.} \\
\thanks{The author would like to thank Osamu Saeki again for private discussions on \cite{saeki} with \cite{kitazawa4}. These discussions continue to encourage the author to continue related studies.

The author was a member of the two projects JSPS KAKENHI Grant Number JP17H06128 and JP22K18267. Principal investigators are both Osamu Saeki. This work was also supported by these projects.
He is also a researcher at Osaka Central
		Advanced Mathematical Institute, supported by MEXT Promotion of Distinctive Joint Research Center Program JPMXP0723833165: he is one of OCAMI researchers whereas he is not employed there. This also helps our studies.

 
All data directly and essentially supporting the present study are in the present paper.}
\section{Preliminary.}
\subsection{Terminologies, notions and notation on smooth or real algebraic manifolds and maps.}
For a (topological) manifold, (a space regarded as a) polyhedron, (one regarded as a) CW complex and (more generally, one regarded as a) cell complex, for example, we can define its dimension uniquely. For such a space $X$, $\dim X$ denotes its dimension, which is a non-negative integer.

For a differentiable manifold $X$, $T_x X$ denotes its tangent space at $x$. 
For a differentiable map $c:X \rightarrow Y$ between differentiable manifolds, $x \in X$ is a {\it singular} point of $c$ if the rank of the differential ${dc}_x:T_x X \rightarrow T_{c(x)} Y$ is smaller than both $\dim X$ and $\dim Y$. 
We define the {\it singular set} $S(c)$ of $c$ as the set consisting of all singular points of $c$.

${\mathbb{R}}^k$ denotes the $k$-dimensional Euclidean space for $k \geq 1$, endowed with the natural differentiable structure and it is a smooth manifold. We can also give the standard metric. 
${\mathbb{R}}^1$ is also denoted by $\mathbb{R}$, which is from natural notation. For $x \in {\mathbb{R}}^k$, $||x|| \geq 0$ denotes the distance between the origin $0 \in {\mathbb{R}}^k$ and $x$ under the metric. 

This is also a $k$-dimensional real algebraic manifold and the {\it $k$-dimensional real affine space}. $S^k:=\{x \in {\mathbb{R}}^{k+1} \mid ||x||=1 \}$ denotes the $k$-dimensional unit sphere for $k \geq 0$. This is a 
smooth compact submanifold (of ${\mathbb{R}}^{k+1}$) and has no boundary. It is connected for $k>0$. 

$S^k$ is also a real algebraic hypersurface of ${\mathbb{R}}^{k+1}$ and of $k$-dimensional.
$D^k:=\{x \mid x \in {\mathbb{R}}^{k}, ||x|| \leq 1\}$ denotes the $k$-dimensional unit disk for $k \geq 1$. This is a $k$-dimensional smooth compact and connected submanifold in ${\mathbb{R}}^k$. 

In our paper, as real algebraic manifolds, we only consider {\it non-singular} real algebraic manifolds. We also consider such a manifold represented as some disjoint union of the zero set of some real polynomial map unless otherwise stated. {\it Non-singular} manifolds are defined naturally via implicit function theorem, applied for the real polynomial maps for the zero sets. Of course the real affine space and the unit sphere are non-singular. Our real algebraic maps are represented as the canonical embeddings into the real affine spaces with canonical projections to some connected components.

\subsection{Graphs and Reeb graphs}
We expect fundamental knowledge including terminologies, notions and notation on graphs. 
In short, a {\it graph} is a $1$-dimensional CW complex the closure of each $1$-cell is homeomorphic to a closed interval. 
An {\it edge} is a $1$-cell of the complex and a {\it vertex} is a $0$-cell of the complex. A graph is also a $1$-dimensional polyhedron.
The set of all edges (vertices) of the complex is the {\it edge} (resp. {\it vertex}) {\it set} of the graph. 
A {\it subgraph} of a graph is a subcomplex of a graph. A subgraph of the graph is also a graph. A {\it finite} graph is a graph whose edge set and vertex set are finite. A connected graph is a graph which is connected as a topological space. 
We only consider finite and connected graphs here.

\begin{Def}
An {\it isomorphism} between two graphs is a (PL or piecewise smooth) homeomorphism between them mapping the vertex set of the graph of the domain onto the vertex set of the graph of the target.
\end{Def}

For a map $c:X \rightarrow Y$ between topological spaces, we can define an equivalence relation on $X$ by the relation that $x_1, x_2 \in X$ are equivalent if and only if $x_1$ and $x_2$ are in a same connected component of some preimage $c^{-1}(y)$. 

\begin{Def}
The quotient space $W_c$ of $X$ defined by the relation is called the {\it Reeb space} of $c$.
\end{Def}

Let $q_c:X \rightarrow W_c$ denote the quotient map. $\bar{c}$ can be defined as a map enjoying the relation $c=\bar{c} \circ q_c$ uniquely. $\bar{c}$ is continuous if $c$ is continuous. This is a fundamental exercise on topological spaces and continuous maps between them. 

\begin{Thm}[\cite{saeki}]
\label{thm:1}
For a smooth function $f:X \rightarrow \mathbb{R}$ on a closed manifold $X$ such that $f(S(f))$ is a finite set, the Reeb space $W_f$ of $f$ is a graph such that a point is a vertex if and only if it contains some singular points of $f$ seen as a connected component of some preimage $f^{-1}(y)$.
\end{Thm}
\begin{Def}
As a generalization, consider a smooth function $f:X \rightarrow \mathbb{R}$ on a manifold $X$ with no boundary, for which we can define such a graph $W_f$.
The graph $W_f$ in the previous theorem is the {\it Reeb graph} of $f$.
\end{Def}


\section{On Main Theorems.}
\subsection{Reviewing explicit construction of real algebraic functions and maps to have desired Reeb graphs.}
We review important arguments from \cite{kitazawa3, kitazawa4, kitazawa5, kitazawa6, kitazawa7}.

For a (finite) set $X$, $|X|$ denotes its size. It is a non-negative integer for a finite set $X$.

For a set $X$ endowed with an order denoted by ${\leq}_X$ and an element $x_0$, let $X_{{\leq}_X x_0}:=\{x \in X \mid x {\leq}_X x_0\}$.
For subsets of $\mathbb{R}$, we consider the natural orders and these orders are denoted by the usual notation $\leq$.

Consider a smooth submanifold $Y$ with no boundary of a smooth manifold $X$.
We can consider the canonical inclusion $T_y Y \subset T_y X$. A {\it normal vector} $v_{{\rm n},y} \in T_y Y \subset T_y X$ at $y \in Y$ is a tangent vector such that for each vector $v_y \in T_y Y \subset T_y X$,  $v_{{\rm n},y} $ is perpendicular to $v_y$. Remember that ${\mathbb{R}}^k$ is endowed with the standard Euclidean metric and that tangent vector spaces there are endowed with the natural inner products.

The following is based on \cite{kitazawa3, kitazawa4, kitazawa5, kitazawa6, kitazawa7}. Especially, it is based on \cite{kitazawa6, kitazawa7}. Note that some are improved and that we omit precise exposition on the improvement. This is not essential in our paper.
\begin{Def}
\label{def:4}
In the $k$-dimensional real affine space ${\mathbb{R}}^k$, let $D$ be an open set there and let $\{S_j\}_{j=1}^l$ be a family of $l>0$ non-singular real algebraic hypersurfaces or real algebraic manifolds of dimension $k-1$ in ${\mathbb{R}}^k$ satisfying the following conditions.
\begin{itemize}
\item $S_j$ is some disjoint union of connected components of the zero set of a real polynomial $f_j$. 
\item For the closure $\overline{D}$ of $D$, $\overline{D}-D \subset {\bigcup}_{j=1}^l S_j$.
\item Consider a subset $A$ consisting of some of the $l$ integers from $1$ to $l$. 
Let $S_{A}:={\bigcap}_{j \in A} S_j$ and $S_{A,0}:=S_{A}-({\bigcup}_{j \in {\mathbb{N}}_{\leq l}- A} S_j)$. $S_{A,0}$ is a smooth submanifold with no boundary. For each point $p \in S_{A,0}$, we can choose a normal vector at $p \in S_j$ for each $j \in A$ and we have a basis of a subspace of $T_p {\mathbb{R}}^k$ whose dimension is $|A|$. Furthermore, we can have the direct sum of this $|A|$-dimensional real vector space and $T_p S_{A,0}$ to have the tangent space $T_p {\mathbb{R}}^k$.  Note that this is on so-called transversality of intersections of the hypersurfaces $S_j$.
\end{itemize}

We say $D$ is said to be a {\it topologically normal convenient domain} or {\it TNC domain}.

Furthermore, if our TNC domain $D$ enjoys the following properties, then $D$ is said to be a {\it normal convenient domain} or {\it NC domain}.
\begin{itemize}
\item $D$ is represented as the intersection of (the closure of) a small open neighborhood $U_D$ of $\overline{D}$ and the set $\{x \in {\mathbb{R}}^k \mid f_j(x)>0\}$.
\item The closure $\overline{D}$ is represented as the intersection of (the closure of) a small open neighborhood $U_D$ of $D$ and the set $\{x \in {\mathbb{R}}^k \mid f_j(x) \geq 0\}$.
\end{itemize}
\end{Def}

For a vertex $v$ of a graph $K$, the {\it degree} of $v$ is the number of edges containing $v$.  
\begin{Thm}[\cite{bodinpopescupampusorea}]
\label{thm:2}
Let $K$ be a finite and connected graph satisfying the following conditions.
\begin{itemize}
\item The degree of each vertex of $K$ is $1$ or $3$.
\item There exists a piecewise smooth function $g_K:K \rightarrow \mathbb{R}$ satisfying the following conditions.
\begin{itemize}
\item On each edge, $g_K$ is an embedding.
\item At distinct vertices, the values of $g_K$ are distinct.
\item If at a vertex $v$, $g_K$ has a local extremum, then $v$ is of degree $1$.
\item We have a piecewise smooth embedding $g_{K,{\mathbb{R}}^2}:K \rightarrow {\mathbb{R}}^2$ such that $g_K$ is the composition of $g_{K,{\mathbb{R}}^2}$ with the projection to the first component of ${\mathbb{R}}^2$. 
\end{itemize}
\end{itemize}
Then we have a TNC domain $D:=D_K$ in ${\mathbb{R}}^2$ enjoying the following properties.
\begin{enumerate}
\item For the closure $\overline{D}$, $\overline{D}-D$ is some disjoint union of non-singular connected real algebraic hypersurfaces $S_j$ of dimension $1$ and has no boundary.  
\item We have a graph $K_D$ enjoying the following properties.
\begin{enumerate}
\item Its underlying space consists of all connected components of preimages considered for the restriction of the projection of ${\mathbb{R}}^2$ to the first component to the subset $\overline{D}$. It is also regarded as the quotient space of $\overline{D}$.
\item A point in the underlying space is a vertex if and only if it is a connected component containing some singular points of the natural smooth function defined as the restriction of the projection of ${\mathbb{R}}^2$ to the first component to the subset $\overline{D}-D$.
\item $K_D$ is isomorphic to $K$. This is given in the following way. First we map each point $p \in K$ by $g_{K,{\mathbb{R}}^2}$. Second we choose a suitable point $p_D \in \overline{D}$ such that the values obtained by mapping this point and the point $g_{K,{\mathbb{R}}^2}(p)$ by using the projection of ${\mathbb{R}}^2$ to the first component are same. We have a point in $K_D$ by mapping $p_D$ by the natural quotient map from $\overline{D}$ onto $K_D$.
\end{enumerate}
\item $D$ is also an NC domain.
\end{enumerate} 
\end{Thm}

We give an example for Definition \ref{def:4} and Theorem \ref{thm:2} fundamental and important in our paper.
\begin{Ex}
\label{ex:1}
Let $t_1<t_2$ be real numbers. FIGURE \ref{fig:1} shows an NC domain in ${\mathbb{R}}^2$ surrounded by a sufficiently large circle and $l_0>0$ circles bounded by mutually disjoint disks whose radii are $\frac{t_2-t_1}{2}$ and which intersect the two straight lines represented as the zero sets of the real polynomials $x_1-t_1$ and $x_1-t_2$. We call such an NC domain an {\it $(l_0,t_1,t_2)$-type domain with circles}.

\begin{figure}
	
	\includegraphics[height=75mm, width=100mm]{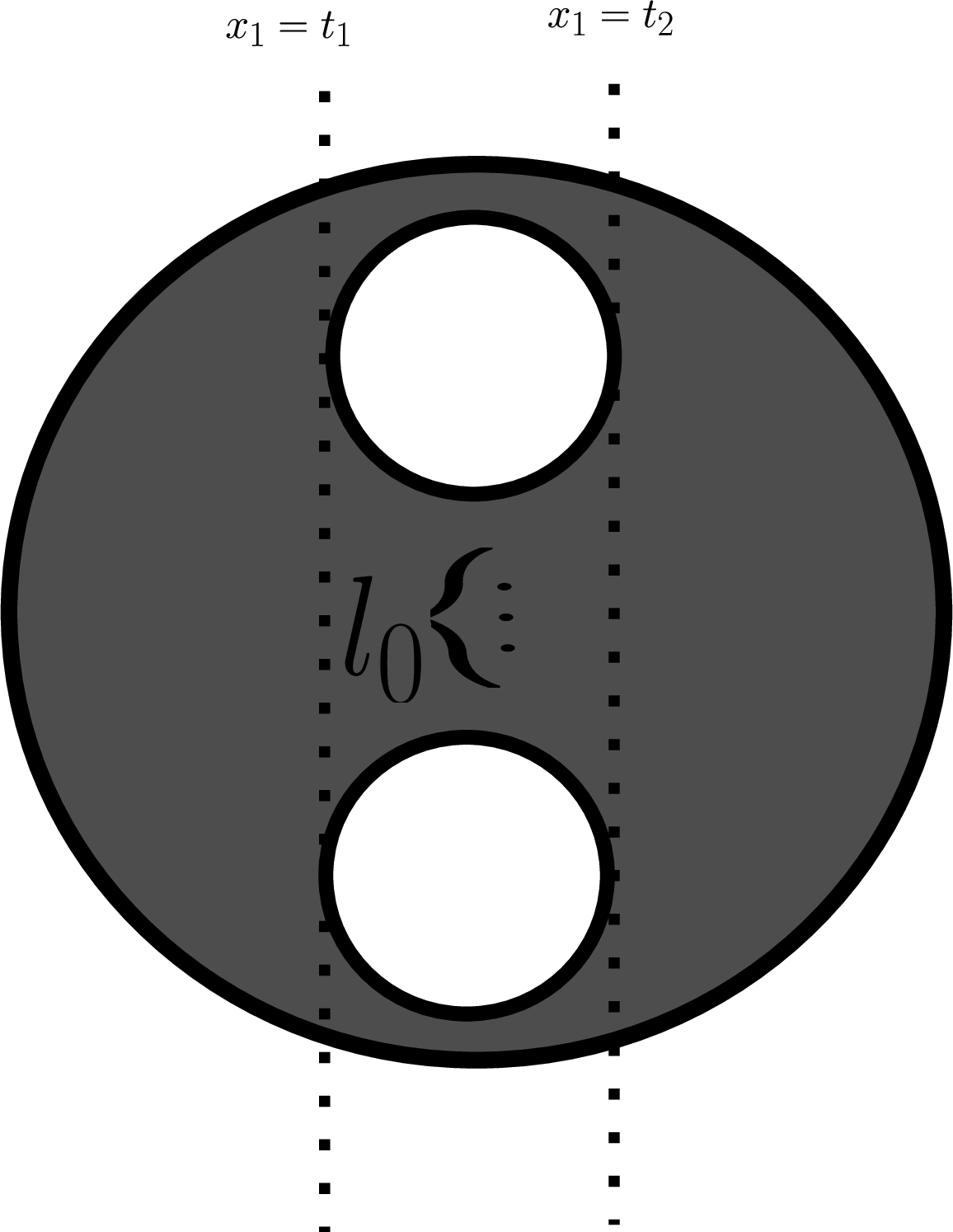}

	\caption{An $(l_0,t_1,t_2)$-type domain with circles.}
	\label{fig:1}
\end{figure}
We can generalize this. Let $\{t_{j,1},t_{j,2}\}_{j=1}^i$ be a sequence of pairs of real numbers of length $i>0$ satisfying the conditions $t_{j,1}<t_{j,2}$ for $1 \leq j \leq i$ and $t_{j^{\prime},2} \leq t_{j^{\prime}+1,1}$ for $1 \leq j^{\prime} \leq i-1$. We can define an NC domain in ${\mathbb{R}}^2$ surrounded by a sufficiently large circle and $l_{j,0}>0$ circles bounded by mutually disjoint disks whose radii are $\frac{t_{j,2}-t_{j,1}}{2}$ and which intersect the straight lines represented as the zero sets of the real polynomials represented as the zero sets of the real polynomials $x_1-t_{j,1}$ and $x_1-t_{j,2}$ for $1 \leq j \leq i$. Here we respect FIGURE 1 of \cite{kitazawa3} partially for example and we have introduced our new definition.

We can consider higher dimensional cases by replacing circles by spheres for example.

Note that for these NC domains $D$ in ${\mathbb{R}}^k$, we can choose the open neighborhood $U_D$ as the Euclidean space ${\mathbb{R}}^k$ for example.
\end{Ex}
\begin{Rem}
For regions here, see also \cite{kohnpieneranestadrydellshapirosinnsoreatelen}. This studies explicit cases in the plane ${\mathbb{R}}^2$. Such regions are also called {\it polypols} there. This is also regarded as a study on some generalized cases of regions in Theorem \ref{thm:2}.
\end{Rem}

\begin{Prop}[\cite{kitazawa6, kitazawa7} (see also \cite{kitazawa3, kitazawa5} for example).]
 \label{prop:1}
We abuse terminologies and notation in Definition \ref{def:4}.
Let $m>0$ be a sufficiently large integer. For an NC domain $D \subset {\mathbb{R}}^k$ in ${\mathbb{R}^k}$,
we have an $m$-dimensional non-singular real algebraic manifold $M_D$, which has no boundary, and a smooth real algebraic map $f_{D}:M_D \rightarrow {\mathbb{R}}^k$ enjoying the following properties.
\begin{enumerate}
\item \label{prop:1.1} The image $f_D(M_D)$ is the closure $\overline{D}$.
\item \label{prop:1.2} For the image of the singular set $S(f_D)$ of $f_D$, $f_D(S(f_D))=\overline{D}-D$.
\item \label{prop:1.3} For each point $p \in {f_D}^{-1}(S_{A,0})$, the image of the differential $d{f_D}_p$ at $p$ and $T_{f_D(p)}S_{A,0}$ agree. 
\item \label{prop:1.4} The preimage of a point in $D$ is diffeomorphic to the product of manifolds diffeomorphic to unit spheres and {\rm (}$m-k${\rm )}-dimensional.
\item \label{prop:1.5} The preimage of a point in $\overline{D}-D$ is a one-point set or a manifold diffeomorphic to the product of manifolds diffeomorphic to unit spheres. In the latter case, the dimension of the preimage is lower than $m-k$.
\item \label{prop:1.6} $M_D$ is connected if $\overline{D}$ is connected. $M_D$ is compact if $\overline{D}$ is compact.
\end{enumerate} 
\end{Prop}
Hereafter, we use the notation of the form $(x_1,\cdots,x_{d_0})$ with $x_i:=(x_{i,1},\cdots,x_{i,d_i})$ for points in $x_i \in {\mathbb{R}}^{d_i}$ and $(x_1,\cdots,x_{d_0}) \in {\mathbb{R}}^{{\Sigma}_{j=1}^{d_0} d_j}$ and coordinates, for example. Of course $d_i$, $d_j$  and $d_0$ here are positive integers. 
\begin{proof}[Reviewing original proofs of Proposition \ref{prop:1}]
We define the set $S_D:=\{(x_1,\cdots x_k,y_1,\cdots y_{l}) \in \overline{D} \times {\mathbb{R}}^{m-k} \subset {\mathbb{R}}^k \times {\mathbb{R}}^{m+l-k}={\mathbb{R}}^{m+l} \mid 1 \leq i \leq l, f_i(x_1,\cdots x_k)-||y_i||^2=0\}$. $y_i$ is regarded as a point in ${\mathbb{R}}^{d_i}$ here with the condition ${\Sigma}_{j=1}^l d_j=m+l-k$ and we can choose the dimensions as positive dimensions suitably since $m$ is sufficiently large. We see that this set is an $m$-dimensional non-singular real algebraic manifold and also represented as some disjoint union of connected components of the zero set of the real polynomial map into ${\mathbb{R}}^l$ obtained canonically from the $l$ real polynomials $f_i(x_1,\cdots x_k)-||y_i||^2$. A main ingredient is implicit function theorem on this polynomial map. We consider several cases. \\
\ \\
Case 1 In the case such that for the point $(x_1,\cdots x_k,y_1,\cdots y_{l}) \in S_D$, $(x_1,\cdots x_k) \in D$.

Since $D$ is an NC domain, for each $f_i(x_1,\cdots x_k)-||y_i||^2$, we consider the partial derivative at the point for each variant $y_{j,j^{\prime}}$ where $j^{\prime}$ is an integer $1 \leq j^{\prime} \leq d_j$. The value is $0$ for $j \neq i$ and it is not zero for some $y_{i,j_i}$ with a suitable integer $1 \leq j_i \leq d_i$. The rank of the map into ${\mathbb{R}}^l$ defined canonically from the real polynomials $f_j(x_1,\cdots x_k)-||y_j||^2$ is $l$ at the point. \\
\ \\
Case 2 In the case such that for the point $(x_1,\cdots x_k,y_1,\cdots y_{l}) \in S_D$, $(x_1,\cdots x_k) \in S_{A,0}$.

Since $D$ is an NC domain, for each $f_i(x_1,\cdots x_k)-||y_i||^2$ with $i \in {\mathbb{N}}_l-A$, we consider the partial derivative at the point for each variant $y_{j,j^{\prime}}$ where $j^{\prime}$ is an integer $1 \leq j^{\prime} \leq d_j$. The value is $0$ for $j \neq i$ and it is not zero for some $y_{i,j_i}$. 

For each $f_i(x_1,\cdots x_k)-||y_i||^2$ with $i \in A$, $y_i$ is the origin. We consider the partial derivative at the point for each variant $y_{j,j^{\prime}}$ where $j^{\prime}$ is an integer $1 \leq j^{\prime} \leq d_j$. The value is always $0$. We consider the map into ${\mathbb{R}}^{|A|}$ defined canonically from the $|A|$ polynomials $f_i(x_1,\cdots x_k)-||y_i||^2$. We consider the partial derivatives by each variant $x_j$ for $1 \leq j \leq k$. We have a matrix of the form $|A| \times k$ consisting of the values of the partial derivatives at the point. By the assumption on the transversality on the intersections of the hypersurfaces $S_j$, the rank is $|A|$.

According to this argument, the rank of the map into ${\mathbb{R}}^l$ defined canonically from the real polynomials $f_j(x_1,\cdots x_k)-||y_j||^2$ is $l$ at the point. \\
\ \\

We consider a point $p_1 \in {\mathbb{R}}^k$ such that $p_1$ is sufficiently close to $D$ and not in the closure $\overline{D}$. 
More precisely, we can formulate this by $p_1 \in U_D-\overline{D}$ for example.
Then according to the assumption, $f_i(p_{1,1},\cdots p_{1,k})<0$ for some $i$, the natural set $S_{D,p_1}:=\{(p_1,p_2) \in {\mathbb{R}}^k \times {\mathbb{R}}^{m+l-k}={\mathbb{R}}^{m+l} \mid 1 \leq i \leq l, f_i(p_{1,1},\cdots p_{1,k})-||p_2||^2=0\}$ is empty. From this with the arguments on Case 1 and Case 2, $S_D$ is regarded as a smooth compact submanifold of dimension $m$ with no boundary in ${\mathbb{R}}^{m+l}$ and some disjoint union of connected components of the zero set of the real polynomial map into ${\mathbb{R}}^l$ defined canonically from the $l$ polynomials.
$f_D$ is defined as the natural projection of $M:=S_D \subset {\mathbb{R}}^{m+l}$ to ${\mathbb{R}}^k$. We have (\ref{prop:1.2}) and (\ref{prop:1.3}) mainly from the arguments on the ranks of the differentials and implicit function theorem before.

We can see (\ref{prop:1.1}),   (\ref{prop:1.4}),  (\ref{prop:1.5}) and  (\ref{prop:1.6}) from the construction easily.

This completes the proof. 
\end{proof}

\subsection{Proving Main Theorems and related remarks.}

We apply some arguments first applied in \cite{kitazawa7}. Related to this, we need notions on "parallel" objects.
First for a Riemannian manifold, we can define the unique natural connection. We can consider mutually parallel tangent vectors and parallel subsets for example. 
For ${\mathbb{R}^k}$, we can define natural mutually independent $k$ tangent vectors at the origin $0$. The $j$-th vector of which is regarded as the vector of ${\mathbb{R}^k}$ whose $j$-th component is $1$ and other components of which are $0$. Let $e_j$ denote this.

Consider ${\mathbb{R}}^k$ or more generally, a product $X:={\prod}_{j \in J} X_j$ indexed by $j \in J$. $J$ is regarded as ${\mathbb{N}}_{\leq k}$ naturally in the case $X:={\mathbb{R}}^k$. Hereafter, for a subset $A \subset J$, let ${\pi}_{X,J,A}$ denote the projection to the components indexed by the elements of the subset $A$.
\begin{proof}[A proof of Main Theorem \ref{mthm:1}]
We define a new open set in ${\mathbb{R}}^k \times {\mathbb{R}}={\mathbb{R}^{k+1}}$. More precisely, we redefine $D$ as a suitable open set there. We abuse the notation in Definition \ref{def:4} for example.

Consider ${\pi}_{{\mathbb{R}}^{k+1},{\mathbb{N}}_{\leq k+1},{\mathbb{N}}_{\leq k}}$ and the preimages of $D$ and the hypersurfaces $S_j$ for this projection.
Consider ${\pi}_{{\mathbb{R}}^{k+1},{\mathbb{N}}_{\leq k+1},\{1,k+1\}}$ and the preimages of an $(i,t_1,t_2)$-type domain with circles and the circles surrounding the domain for this projection.
We can redefine our new open set $D$ as the intersection of the two preimages of the two NC domains. We can also redefine our new hypersurfaces surrounding the new domain as the preimages of the original hypersurfaces. Let the new hypersurface denote $S_{j}$ for the preimage of $S_j$ in the original situation and $S_{l+j}$ for the preimage of the $j$-th circle in the family of the $i$-circles surrounding the $(i,t_1,t_2)$-type domain with circles. The outermost sufficiently large circle is denoted by $S_{l+i+1}$.

We explain about the fact that our new set $D$ is an NC domain. We can easily check the conditions on the real polynomials and the zero sets and the fact that the hypersurfaces are non-singular, by our construction. The preimages are regarded as the natural cylinders for the original hypersurfaces. We check the condition on the transversality of the intersections of hypersurfaces.

$S_{j_1}$ and $S_{j_2}$ do not intersect for $j_1<j_2 \leq l$ or $l+1 \leq j_1<j_2 \leq l+(i+1)=l+i+1$.

We consider the remaining case. Suppose that $S_{j_1} \bigcap S_{j_2}$ is not empty for some $j_1<j_2$ and that $p$ is a point in this intersection. Let $p_1$ denote the value obtained by mapping this by the projection to the first component. 
We can consider the following exactly three cases by the assumption or the condition C4. \\
\ \\
Case A The preimage ${g_{K_D}}^{-1}(p_1)$ for the original function $g_{K_D}$ contains no vertices and $p_1 \neq t_j$ for $j=1,2$.  \\
A normal vector at $p \in S_{j_1}$ must be parallel to a vector of the form ${\Sigma}_{j=1}^k t_je_j$ for some $j \neq 1$ and $t_j \neq 0$.
A normal vector at $p \in S_{j_2}$ must be parallel to a vector of the form $t_1e_1+t_{k+1}e_{k+1}$ with $t_{k+1} \neq 0$. Remember that the outermost circle in the $(i,t_1,t_2)$-type domain with circles is chosen sufficiently large and that $S_{l+i+1}$ is also a sufficiently large cylinder of the circle. \\
 \ \\
Case B The preimage ${g_{K_D}}^{-1}(p_1)$ for the original function $g_{K_D}$ contains no vertices and $p_1=t_j$ for either $j=1$ or $j=2$.  \\
A normal vector at $p \in S_{j_1}$ must be parallel to a vector of the form ${\Sigma}_{j=1}^k t_j e_j$ for some $j \neq 1$ and $t_j \neq 0$.
A normal vector at $p \in S_{j_2}$ must be parallel to a vector of the form $t_1e_1$ for some $t_1 \neq 0$. \\
\ \\
Case C The preimage ${g_{K_D}}^{-1}(p_1)$ for the original function $g_{K_D}$ contains some vertices and $p_1 \neq t_j$ for $j=1,2$. \\
In the case $p$ is not a singular point of the function of the restriction of the projection of ${\mathbb{R}}^k$ to the first component to "the original set $\overline{D}-D$", a normal vector at $p \in S_{j_1}$ must be parallel to a vector of the form ${\Sigma}_{j=1}^k t_j e_j$ for some $j \neq 1$ and $t_j \neq 0$.
$p$ may be a singular point of the function of the restriction of the projection of ${\mathbb{R}}^k$ to the first component to "the original set $\overline{D}-D$" and in this case, a normal vector at $p \in S_{j_1}$ must be parallel to a vector of the form $t_1e_1$ with $t_1 \neq 0$.
A normal vector at $p \in S_{j_2}$ must be parallel to a vector of the form $t_1e_1+t_{k+1}e_{k+1}$ with $t_{k+1} \neq 0$. Remember again that the outermost circle in the $(i,t_1,t_2)$-type domain with circles is chosen sufficiently large and that $S_{l+i+1}$ is also a sufficiently large cylinder of the circle. \\\\
\ \\
We can see the transversality easily from this. We also see that our new $D$ is an NC domain. $U_D$ may be large. However it is regarded as a desired open neighborhood of our new $D$. 

We apply Proposition \ref{prop:1} to have a new map $f_D:M_D \rightarrow {\mathbb{R}}^{k+1}$ whose image is the closure $\overline{D}$ of the new set $D$.
This also gives a part of our proof of Main Theorem \ref{mthm:2} and a part of our proof of Main Theorem \ref{mthm:3} essentially.

We consider the composition of $f_D$ with the projection to the first component and its Reeb graph. We have a smooth real algebraic function, denoted by $f_i:M_i \rightarrow {\mathbb{R}}^{k+1}$ with $\dim M_i=m_i$. Let us see that $f_i$ is our desired function. First, we have the following by fundamental arguments.
\begin{itemize}
\item Remember that $\bar{f_i}:W_{f_i} \rightarrow \mathbb{R}$ enjoys $f_i=\bar{f_i} \circ q_{f_i}$. $\bar{f_i}$ is represented as the composition of a suitable piecewise smooth map ${\phi}_{f_i}$ onto $K_D$ with the function $g_{K_D}$.
\item The restriction of ${\phi}_{f_i}$ to the preimage ${{\phi}_{f_i}}^{-1}({g_{K_D}}^{-1}(\mathbb{R}-(t_1,t_2)))$ of the complementary set of the open interval $(t_1,t_2)$ in $\mathbb{R}$ is regarded to give an isomorphism of the two graphs.
\item The restriction of ${\phi}_{f_i}$ to the preimage ${{\phi}_{f_i}}^{-1}({g_{K_D}}^{-1}((t_1,t_2)))$ of the open interval $(t_1,t_2)$ is regarded as an ($i+1$)-fold covering.
\item Consider the set of all vertices of the graph $W_{f_i}$ in the preimage ${{\phi}_{f_i}}^{-1}({g_{K_D}}^{-1}((t_1,t_2)))$ and compare this to the set of all points in the preimages of some vertices in ${g_{k_D}}^{-1}((t_1,t_2))$. The size of this set is $i+1$ times the size of the number of all vertices in ${g_{k_D}}^{-1}((t_1,t_2))$.  
\end{itemize}
We can easily see these properties. We explain about the third and the fourth properties here mainly. 
By the symmetry, it is sufficient to consider the case C5.1. The preimage of $\{p_{t_1,1},p_{t_1,2}\}$ considered for ${\phi}_{f_i}$ is a two-point set with the discrete topology and consists of two vertices. The preimage of $\{p_{t_2}\}$ considered for ${\phi}_{f_i}$ is a one-point set and consists of some vertex. Consider the preimage of the union of the two given embedded arcs $e_{t_1,t_2,1}$ and $e_{t_1,t_2,2}$. The restriction of the map ${\phi}_{f_i}$ to the interior of this 1-dimensional set is regarded as an ($i+1$)-fold covering. The interior of this $1$-dimensional set is obtained by removing the preimage of $\{p_{t_1,1},p_{t_1,2}\}$ considered for ${\phi}_{f_i}$ and the preimage of $\{p_{t_2}\}$ considered for ${\phi}_{f_i}$.

We consider the case C5.1 further to complete the proof. Consider a (presented) small connected embedded arc in $e_{t_1,t_2,1} \bigcap e_{t_1,t_2,2}$ containing $p_{t_2}$ and remove the point $p_{t_2}$. We consider the preimage of the resulting arc for ${\phi}_{f_i}$. There ${\phi}_{f_i}$ is an ($i+1$)-fold covering.
Assume that $\bar{f_i}:W_{f_i} \rightarrow \mathbb{R}$ is represented as the composition of some embedding $g_{K_{D,t_1,t_2,i},{\mathbb{R}}^2}:W_{f_i} \rightarrow {\mathbb{R}}^2$ with the projection to the first component and we explain about the contradiction. We have $i+1$ distinct points in the image $g_{K_{D,t_1,t_2,i},{\mathbb{R}}^2}({\rm Int}\ (e_{t_1,t_2,1} \bigcap e_{t_1,t_2,2})) \subset {\mathbb{R}}^2$ of the interior of the intersection $e_{t_1,t_2,1} \bigcap e_{t_1,t_2,2}$ where the values of the projection to the first component are a same real number $t_{1,2} \in (t_1,t_2)$. Let $p_{e,j}$ denote the $j$-th point of these $i+1$ points. We also have $i+1$ pairs of embedded arcs in $W_{f_i}$ with the following.
\begin{itemize}
\item For the $j$-th pair, one arc connects $p_{t_1,1}$ and $p_{e,j}$ and the other arc connects $p_{t_1,2}$ and $p_{e,j}$.
\item As presented, the values of the projection to the first component at $p_{e,j} \in {\mathbb{R}}^2$ are all $t_{1,2}$. 
\end{itemize}
We can see the contradiction, implying that $g_{K_{D,t_1,t_2,i_1},{\mathbb{R}}^2}$ is not an embedding.

This completes the proof.
\end{proof}
\begin{MainThm}
\label{mthm:2}
Let $D \subset {\mathbb{R}}^k$ be an open and connected set satisfying the following conditions. We abuse the notation in Definition \ref{def:4} for example. We also respect the conditions in Main Theorem \ref{mthm:1}.
\begin{itemize}
\item[{C1}$^{\prime}$] The closure $\overline{D}$ is compact and connected. $\overline{D}-D$ is some disjoint union of finitely many {\it non-singular} real algebraic hypersurfaces $S_j$ of dimension $k-1$, indexed by $j \in J$ in a finite set $J$ of size $l>0$, and has no boundary.
\item[{C2}$^{\prime}$] $D$ is, as in Definition \ref{def:4}, an {\it NC domain}.
\item[{C3}$^{\prime}$] There exists a finite and connected graph $K_D$ enjoying the following properties.
\begin{itemize}
\item[{C3.1}$^{\prime}$] Its underlying space consists of all connected components of preimages considered for the restriction of the projection of ${\mathbb{R}}^k$ to the first component to the subset $\overline{D}$. It is also regarded as the quotient space of $\overline{D}$.
We also have a natural map $g_{K_D}:K_D \rightarrow \mathbb{R}$ by considering the value of the restriction of the projection of ${\mathbb{R}}^k$ to the first component at each point in $p \in K_D$, representing some preimage.
\item[{C3.2}$^{\prime}$] A point in the underlying space is a vertex if and only if it is a connected component containing some singular points of the natural smooth function defined as the restriction of the projection of ${\mathbb{R}}^k$ to the first component to the subset $\overline{D}-D$.
\item[{C3.3}$^{\prime}$] $g_{K_D}$ is a piecewise smooth function and injective on each edge of the graph $K_D$.
\end{itemize}
\end{itemize}
Suppose that two distinct values $t_1<t_2$ satisfy the following conditions.
\begin{itemize}
\item[{C4}$^{\prime}$] ${g_{K_D}}^{-1}(t_1)$ and ${g_{K_D}}^{-1}(t_2)$ contain no vertices of $K_D$.
\item[{C5}$^{\prime}$] ${g_{K_D}}^{-1}(t_1)$ contains three distinct points ${p_{t_1,1}}^{\prime}$, ${p_{t_1,2}}^{\prime}$ and ${p_{t_1,3}}^{\prime}$. ${g_{K_D}}^{-1}(t_2)$ contains three distinct points ${p_{t_2,1}}^{\prime}$, ${p_{t_2,2}}^{\prime}$ and ${p_{t_2,3}}^{\prime}$. 
We have some embedded arc $e_{t_1,t_2,j_1,j_2}$ in the graph $K_D$ connecting ${p_{t_1,j_1}}^{\prime}$ and ${p_{t_2,j_2}}^{\prime}$ making the image $g_{K,D}({\rm Int}\ e_{t_1,t_2,j_1,j_2})$ of the interior of each arc and the open interval $(t_1,t_2)$ agree for $1 \leq j_1,j_2 \leq 3$. 
\end{itemize}
Then we have a family $\{{K_{D,t_1,t_2,i}}^{\prime}\}$ of graphs indexed by positive integers $i$ enjoying the following properties.
\begin{enumerate}
\item ${K_{D,t_1,t_2,i_1}}^{\prime}$ and ${K_{D,t_1,t_2,i_2}}^{\prime}$ are not isomorphic for $i_1 \neq i_2$.
\item For a sufficiently large integer ${m_i}^{\prime}>0$, we have a suitable ${m_i}^{\prime}$-dimensional non-singular real algebraic closed and connected manifold ${M_i}^{\prime}$ and a smooth real algebraic function ${f_i}^{\prime}:{M_i}^{\prime} \rightarrow \mathbb{R}$ whose Reeb graph $W_{{f_i}^{\prime}}$ is isomorphic to ${K_{D,t_1,t_2,i}}^{\prime}$.
\item We cannot embed ${K_{D,t_1,t_2,i}}^{\prime}$ into ${\mathbb{R}}^2$.
\end{enumerate}
\end{MainThm}
\begin{proof}
We abuse the notation as in our proof of Main Theorem \ref{mthm:1}.
As presented in our proof of Main Theorem \ref{mthm:1}, we can similarly obtain a smooth real algebraic map $f_D:M_D \rightarrow {\mathbb{R}}^{k+1}$ by using Proposition \ref{prop:1}. We consider the composition of $f_D$ with the projection to the first component and its Reeb graph. We have a smooth real algebraic function, denoted by $f_i:{M_i}^{\prime} \rightarrow {\mathbb{R}}^{k+1}$ with $\dim {M_i}^{\prime}={m_i}^{\prime}$.
We should note that instead of using an $(i,t_1,t_2)$-type domain with circles, we use an $(i+7,t_1,t_2)$-type domain for obtaining our new open set $D$ and map $f_D$.

We respect the condition {C4}$^{\prime}$ to have the following. 
\begin{itemize}
\item Remember that $\bar{{f_i}^{\prime}}:W_{{f_i}^{\prime}} \rightarrow \mathbb{R}$ enjoys ${f_i}^{\prime}=\bar{{f_i}^{\prime}} \circ q_{{f_i}^{\prime}}$. $\bar{{f_i}^{\prime}}$ is represented as the composition of a suitable piecewise smooth map ${\phi}_{{f_i}^{\prime}}$ onto $K_D$ with the function $g_{K_D}$.
\item The restriction of ${\phi}_{{f_i}^{\prime}}$ to the preimage ${{\phi}_{{f_i}^{\prime}}}^{-1}({g_{K_D}}^{-1}(\mathbb{R}-(t_1,t_2)))$ of the complementary set of the open interval $(t_1,t_2)$ in $\mathbb{R}$ is regarded to give an isomorphism of the two graphs.
\item The restriction of ${\phi}_{{f_i}^{\prime}}$ to the preimage ${{\phi}_{{f_i}^{\prime}}}^{-1}({g_{K_D}}^{-1}((t_1,t_2)))$ of the open interval $(t_1,t_2)$ is regarded as an ($i+8$)-fold covering.
\item Consider the set of all vertices of the graph $W_{{f_i}^{\prime}}$ in the preimage ${{\phi}_{{f_i}^{\prime}}}^{-1}({g_{K_D}}^{-1}((t_1,t_2)))$ and compare this to the set of all points in the preimages of some vertices in ${g_{k_D}}^{-1}((t_1,t_2))$. The size of this set is $i+8$ times the size of the number of all vertices in ${g_{k_D}}^{-1}((t_1,t_2))$.  
\end{itemize}
We can easily see these properties. We explain about the third and the fourth properties here mainly. We respect the condition {C5}$^{\prime}$. 
The preimage of $\{{p_{t_1,1}}^{\prime},{p_{t_1,2}^{\prime}},{p_{t_1,3}}^{\prime}\}$ considered for ${\phi}_{{f_i}^{\prime}}$ is a three-point set with the discrete topology and consists of three vertices. The preimage of $\{{p_{t_2,1}}^{\prime},{p_{t_2,2}^{\prime}},{p_{t_2,3}}^{\prime}\}$ considered for ${\phi}_{{f_i}^{\prime}}$ is a three-point set with the discrete topology and consists of three vertices. Consider the preimage of the union of the given embedded arcs $e_{t_1,t_2,j_1,j_2}$. The restriction of the map ${\phi}_{{f_i}^{\prime}}$ to the interior of this 1-dimensional set is regarded as an ($i+8$)-fold covering. The interior of this $1$-dimensional set is obtained by removing the preimage of $\{{p_{t_1,1}}^{\prime},{p_{t_1,2}}^{\prime},{p_{t_1,3}}^{\prime}\}$ considered for ${\phi}_{{f_i}^{\prime}}$ and the preimage of $\{{p_{t_2,1}}^{\prime},{p_{t_2,2}}^{\prime},{p_{t_2,3}^{\prime}}\}$ considered for ${\phi}_{{f_i}^{\prime}}$. From this argument, we can find a so-called $K_{3,3}$-graph as a subgraph in the graph ${K_{D,t_1,t_2,i}}^{\prime}$. It is well-known that we cannot embed such a graph into ${\mathbb{R}}^2$.
This completes the proof.
\end{proof}
\begin{MainThm}
\label{mthm:3}
Let $D \subset {\mathbb{R}}^k$ be an open and connected set satisfying the following conditions. We abuse the notation in Definition \ref{def:4} for example and respect some conditions in Main Theorems \ref{mthm:1} and \ref{mthm:2}.
\begin{itemize}
\item[{C1}$^{\prime \prime}$] The conditions {C1}$^{\prime}$, {C2}$^{\prime}$ and {C3}$^{\prime}$ hold.
\item[{C2}$^{\prime \prime}$] There exist two distinct real numbers $t_1<t_2$ and ${g_{K_D}}^{-1}(t_j)$ contains no vertices of $K_D$ for $j=1,2$.
\item[{C3}$^{\prime \prime}$] ${g_{K_D}}^{-1}(t_j)$ contains a subset $\{{p_{t_j,j^{\prime}}}^{\prime \prime}\}_{j^{\prime}=1}^{a_j}$ consisting of exactly $a_j$ distinct points where $a_1=2$ and $a_2=3$. We have some embedded arcs in the graph $K_D$ in the following.
\begin{itemize}
\item[{C3.1}$^{\prime \prime}$] Some embedded arc ${e_{t_1,t_2,j_1,j_2}}^{\prime}$ in the graph $K_D$ connecting ${p_{t_1,j_1}}^{\prime \prime}$ and ${p_{t_2,j_2}}^{\prime \prime}$ making the image $g_{K_D}({\rm Int}\ {e_{t_1,t_2,j_1,j_2}}^{\prime})$ of the interior of the arc and the open interval $(t_1,t_2)$ agree for each $j_1$ and $j_2$. 
\item[{C3.2}$^{\prime \prime}$] Some embedded arc ${e_{t_1,t_1,1,2}}^{\prime}$ in the graph $K_D$ connecting ${p_{t_1,1}}^{\prime \prime}$ and ${p_{t_1,2}}^{\prime \prime}$ putting the image $g_{K_D}({\rm Int}\ {e_{t_1,t_1,1,2}}^{\prime})$ of the interior of the arc outside the subset $\{t_1,t_2\}$.
\item[{C3.2}$^{\prime \prime}$] Some embedded arc ${e_{t_2,t_2,j_1,j_2}}^{\prime}$ in the graph $K_D$ connecting ${p_{t_2,j_1}}^{\prime \prime}$ and ${p_{t_2,j_2}}^{\prime \prime}$ putting the image $g_{K_D}({\rm Int}\ {e_{t_2,t_2,j_1,j_2}}^{\prime})$ of the interior of the arc outside the subset $\{t_1,t_2\}$ for each $j_1<j_2$.
\end{itemize}
\end{itemize}
Then we have a family $\{K_{D,t_1,t_2,i_1,i_2,i_3}\}$ of graphs indexed by positive integers $i_1$, $i_2$ and $i_3$ enjoying the following properties.
\begin{enumerate}
\item $K_{D,t_1,t_2,i_{1,1},i_2,i_3}$ and $K_{D,t_1,t_2,i_{1,2},i_2.i_3}$ are not isomorphic for $i_{1,1} \neq i_{1,2}$.
\item $K_{D,t_1,t_2,i_1,i_{2,1},i_3}$ and $K_{D,t_1,t_2,i_1,i_{2,2}.i_3}$ are not isomorphic for $i_{2,1} \neq i_{2,2}$.
\item $K_{D,t_1,t_2,i_1,i_2,i_{3,1}}$ and $K_{D,t_1,t_2,i_1,i_2.i_{3,2}}$ are not isomorphic for $i_{3,1} \neq i_{3,2}$.
\item For a sufficiently large integer $m_{i_1,i_2,i_3}>0$, we have a suitable $m_{i_1,i_2,i_3}$-dimensional non-singular real algebraic closed and connected manifold $M_{i_1,i_2,i_3}$ and a smooth real algebraic function $f_{i_1,i_2,i_3}:M_{i_1,i_2,i_3} \rightarrow \mathbb{R}$ whose Reeb graph $W_{f_{i_1,i_2,i_3}}$ is isomorphic to $K_{D,t_1,t_2,i_1,i_2,i_3}$.
\item We cannot embed each graph $K_{D,t_1,t_2,i_1,i_2,i_3}$ into ${\mathbb{R}}^2$.
\end{enumerate}
\end{MainThm}
\begin{proof}
	We abuse the notation in our proofs of Main Theorems \ref{mthm:1} and \ref{mthm:2} for example. Our main story of our proof is similar to those of Main Theorems \ref{mthm:2} and \ref{mthm:3}.
	
	As the domain to obtain our new open set $D$ and map $f_D$, we use an domain for $\{t_{j,1},t_{j,2}\}_{j=1}^3$ in Example \ref{ex:1} instead in the following way.
\begin{itemize}
\item $t_{1,1}$ is chosen as a sufficiently small real number ${t_1}^{\prime}<t_1$.
\item $t_{1,2}=t_{2,1}=t_1$.
\item $t_{2,2}=t_{3,1}=t_2$.
\item $t_{3,2}$ is chosen as a sufficiently large real number ${t_2}^{\prime}>t_2$.
\item $l_{1,0}=i_1$.
\item $l_{2,0}=i_2+8$.
\item $l_{3,0}=i_3+1$
\end{itemize}

We respect the condition {C3}$^{\prime \prime}$ to have the following. 
\begin{itemize}
\item Remember that $\bar{f_{i_1,i_2,i_3}}:W_{{f_i}_{1,2,3}} \rightarrow \mathbb{R}$ enjoys $f_{i_1,i_2,i_3}=\bar{f_{i_1,i_2,i_3}} \circ q_{f_{i_1,i_2,i_3}}$. $\bar{f_{i_1,i_2,i_3}}$ is represented as the composition of a suitable piecewise smooth map ${\phi}_{f_{i_1,i_2,i_3}}$ onto $K_D$ with the function $g_{K_D}$.
\item The restriction of ${\phi}_{f_{i_1,i_2,i_3}}$ to the preimage ${{\phi}_{f_{i_1,i_2,i_3}}}^{-1}({g_{K_D}}^{-1}((t_1,t_2)))$ of the open interval $(t_1,t_2)$ is regarded as an ($l_{2,0}+1$)-fold covering.
The restriction of ${\phi}_{f_{i_1,i_2,i_3}}$ to the preimage ${{\phi}_{f_{i_1,i_2,i_3}}}^{-1}({g_{K_D}}^{-1}([t_0,t_1)))$ of the interval $[t_0,t_1)$ is regarded as an ($l_{1,0}+1$)-fold covering where a sufficiently small number $t_0$ is chosen suitably satisfying ${t_1}^{\prime}<t_0<t_1$.
The restriction of ${\phi}_{f_{i_1,i_2,i_3}}$ to the preimage ${{\phi}_{f_{i_1,i_2,i_3}}}^{-1}({g_{K_D}}^{-1}((t_2,t_3]))$ of the interval $(t_2,t_3]$ is regarded as an ($l_{3,0}+1$)-fold covering where a sufficiently large number $t_3$ is chosen suitably satisfying $t_2<t_3<{t_2}^{\prime}$.
\item 
Let $A_1:=[t_0,t_1)$, $A_2:=(t_1,t_2)$ and $A_3:=(t_2,t_3]$. Consider the set of all vertices of the graph $W_{f_{i_1,i_2,i_3}}$ in the preimage ${{\phi}_{f_{i_1,i_2,i_3}}}^{-1}({g_{K_D}}^{-1}(A_j))$ and compare this to the set of all points in the preimages of some vertices in ${g_{k_D}}^{-1}(A_j)$ for each $j=1,2,3$. The size of this set is $l_{j,0}+1$ times the size of the number of all vertices in ${g_{k_D}}^{-1}(A_j)$.  
\end{itemize}
As in Main Theorems \ref{mthm:1} and \ref{mthm:2}, we can have smooth real algebraic functions and their Reeb graphs. By considering the several explicitly given embedded arcs as in our proofs of Main Theorems \ref{mthm:1} and \ref{mthm:2}, so-called $K_5$-graphs can be found as subgraphs of the resulting Reeb graphs. It is well-known that we cannot embed $K_5$-graphs into ${\mathbb{R}}^2$. This completes the proof.
\end{proof}
\begin{Rem}
\label{rem:1}
In the situation of Theorem \ref{thm:2}, we can apply Proposition \ref{prop:1} and have a similar map $f_D:M \rightarrow {\mathbb{R}}^k$ similarly. We compose the resulting map with the projection to the first component similarly. As a result, we can apply an argument like one in Case C in our proof of Main Theorem \ref{mthm:2} have a smooth real algebraic function whose Reeb graph is isomorphic to "$K$ and $K_D$ in Theorem \ref{thm:2}". 

Furthermore, for the resulting function "$f$ in Theorem \ref{thm:2}", the function $\bar{f}:W_f \rightarrow \mathbb{R}$, enjoying the relation $f=\bar{f} \circ q_f$, is naturally represented as the composition of some piecewise smooth embedding into ${\mathbb{R}}^2$ with the projection to the first component. "Resulting functions in Main Theorems" are essentially different from this function with respect to the fact that "$\bar{f}:W_f \rightarrow \mathbb{R}$ in Main Theorems" cannot be represented as the compositions of embeddings into ${\mathbb{R}}^2$ with the projection to the first component.
\end{Rem}
We present domains and graphs related to Main Theorems.
\begin{Ex}
\label{ex:2}
The first figure in FIGURE \ref{fig:2} is a graph obtained by considering a $(1,t_{0,1},t_{0,2})$-type domain with circles, denoted by $D$ here, and the function $g_{K,D}$ which is naturally obtained, another $(1,t_{1},t_{2})$-type domain with circles satisfying $t_{0,1}<t_1<t_{0,2}<t_2$, and Main Theorem \ref{mthm:1}. The edge represented by the left gray arc goes over the edge represented by the separated two arcs. This graph also respects the function $g_{K_D}$ on the graph $K_D$. The second graph shows a graph isomorphic to the first graph and this is embedded into the plane. This graph is regarded as a graph the third graph collapses to. More precisely, by eliminating two edges, we have the second graph. The third graph $K$ respects a natural piecewise smooth function $g_K$ for Theorem \ref{thm:2}.
\begin{figure}
	
	\includegraphics[height=75mm, width=100mm]{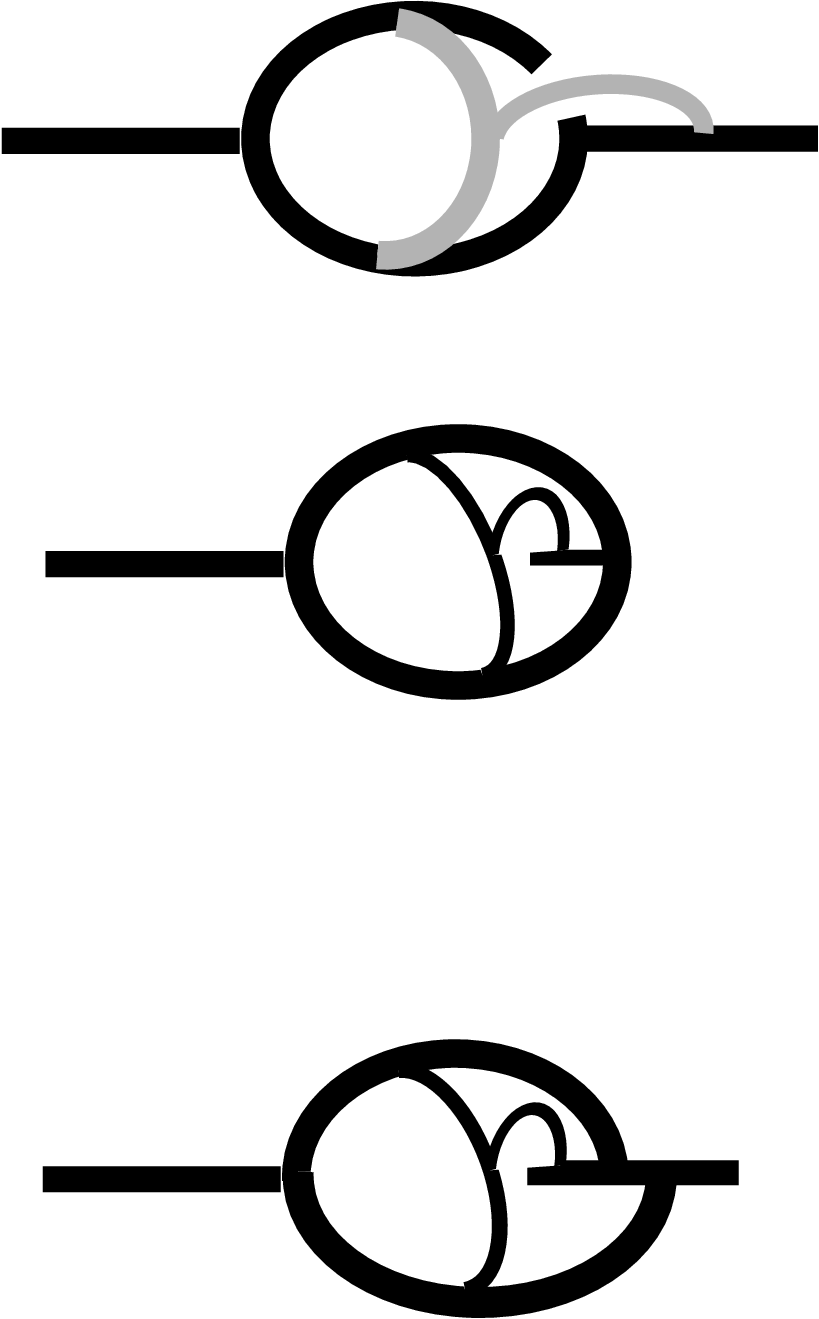}

	\caption{Graphs related to Main Theorem \ref{mthm:1}.}
	\label{fig:2}
\end{figure}
\end{Ex}
Related to Example \ref{ex:2}, we present our natural problem.
\begin{Prob}
Consider a finite and connected graph $K_0$ with at least one edge.
What is the simplest or a most natural smooth real algebraic function on a non-singular real algebraic manifold whose Reeb graph is isomorphic to a graph $K$ collapsing to $K_0$?
\end{Prob}
We can also consider similar variants according to our situations.

For example, it is natural to consider that for a graph $K_0$ with exactly one edge and two vertices, the restriction of the projection of $\mathbb{R}^{k}$ to the first component to the unit sphere $S^{k-1}$ is a desired function where $k \geq 2$.

\begin{Ex}
\label{ex:3}

We can consider a domain in Example \ref{ex:1}, considered in the case $\{t_{j,1},t_{j,2}\}_{j=1}^2$ with $t_{1,1}<t_1<t_{1,2} \leq t_{2,1}<t_2<t_{2,2}$ and $l_{1,0}=l_{2,0}=2$ as the originally given NC domain to apply Main Theorem \ref{mthm:2}.

\end{Ex}

\begin{Rem}
	\label{rem:2}
We give short remarks on the proof of Main Theorem \ref{mthm:3} by presenting explicit cases.  
If the image $g_{K_D}({\rm Int}\ {e_{t_1,t_1,1,2}}^{\prime})$ of the interior of the arc is in the subset $({t_1}^{\prime},t_1)$, then $l_{2,0}=i_2+8$ can be replaced by $l_{2,0}=i_2+6$. If the image $g_{K_D}({\rm Int}\ {e_{t_2,t_2,j_1,j_2}}^{\prime})$ of the interior of the arc is in the subset $(t_2,{t_2}^{\prime})$ for each $j_1<j_2$, then $l_{2,0}=i_2+8$ can be replaced by $l_{2,0}=i_2+5$. If both of the assumptions here are satisfied, then $l_{2,0}=i_2+8$ can be replaced by $l_{2,0}=i_2+4$.

\end{Rem}
\begin{Ex}
	\label{ex:4}
	
	We can consider a domain in Example \ref{ex:1}, considered in the case $\{t_{j,1},t_{j,2}\}_{j=1}^2$ with $t_{1,1}<t_1<t_{1,2} \leq t_{2,1}<t_2<t_{2,2}$ and $(l_{1,0},l_{2,0})=(1,2)$ as the originally given NC domain to apply Main Theorem \ref{mthm:2}. Such a case is regarded as examples for all cases in Remark \ref{rem:2}.
	
\end{Ex}
\begin{Rem}
\label{rem:3}
On NC domains $D$ in Main Theorems, we can consider general NC domains as in Definition \ref{def:4} and discuss our generalized cases. We omit arguments. We can also have various cases only considering (very explicit) cases of Theorem \ref{thm:2} such as cases of Example \ref{ex:1} and higher dimensional versions of the cases.
\end{Rem}


\begin{thebibliography}{25} 
\bibitem{bodinpopescupampusorea} A. Bodin, P. Popescu-Pampu and M. S. Sorea, \textsl{Poincar\'e-Reeb graphs of real algebraic domains}, Revista Matem\'atica Complutense, https://link.springer.com/article/10.1007/s13163-023-00469-y, arXiv:2207.06871v2.
\bibitem{kitazawa1} N. Kitazawa, \textsl{On Reeb graphs induced from smooth functions on $3$-dimensional closed orientable manifolds with finitely many singular values}, Topological Methods in Nonlinear Analysis Vol. 59 No. 2B (2022), 897--912, https://doi.org/10.12775/TMNA.2021.044, arXiv:1902:8841.
\bibitem{kitazawa2} N. Kitazawa, \textsl{On Reeb graphs induced from smooth functions on closed or open surfaces}, Methods of Functional Analysis and Topology Vol. 28 No. 2 (2022), 127--143, doi.org/10.31392/MFAT-npu26\_2.2022.05, arXiv:1908.04340.
\bibitem{kitazawa3} N. Kitazawa, \textsl{Real algebraic functions on closed manifolds whose Reeb graphs are given graphs}, a positive report for publication has been announced to have been sent and this will be published in Methods of Functional Analysis and Topology, arXiv:2302.02339.

\bibitem{kitazawa4} N. Kitazawa, \textsl{Explicit construction of explicit real algebraic functions and real algebraic manifolds via Reeb graphs}, Algebraic and geometric methods of analysis 2023 “The book of abstracts”, 49—51, this is the abstract book of the conference "Algebraic and geometric methods of analysis 2023" (https://www.imath.kiev.ua/$\sim$topology/conf/agma2023/), https://imath.kiev.ua/$\sim$topology/conf/agma2023/contents/abstracts/texts/kitazawa/kitazawa.pdf.

\bibitem{kitazawa5} N. Kitazawa, \textsl{Construction of real algebraic functions with prescribed preimages}, submitted to a refereed journal, arXiv:2303.00953.
\bibitem{kitazawa6} N. Kitazawa, \textsl{Natural real algebraic maps of non-positive codimensions with prescribed images whose boundaries consist of non-singular real algebraic hypersurfaces intersecting with transversality}, submitted to a refereed journal, arXiv:2303.10723.
\bibitem{kitazawa7} N. Kitazawa, \textsl{Explicit smooth real algebraic functions which may have both compact and non-compact preimages on smooth real algebraic manifolds}, arXiv:2304.02372, this version is submitted to a refereed journal.

\bibitem{kobayashisaeki} M. Kobayashi and O. Saeki, \textsl{Simplifying stable mappings into the plane from a global viewpoint}, Trans. Amer. Math. Soc. 348 (1996), 2607--2636.
\bibitem{kohnpieneranestadrydellshapirosinnsoreatelen} K. Kohn, R. Piene, K. Ranestad, F. Rydell, B. Shapiro, R. Sinn, M-S. Sorea and S. Telen, \textsl{Adjoints and Canonical Forms of Polypols}, arXiv:2108.11747.
\bibitem{masumotosaeki} Y. Masumoto and O. Saeki, \textsl{A smooth function on a manifold with given Reeb graph}, Kyushu J. Math. 65 (2011), 75--84; doi: 10.2206/kyushujm.65.75. 
\bibitem{matsuzaki} S. Matsuzaki, \textsl{Multiplicity of Finite Graphs Over the Real Line}, Tokyo J. Math. 37 (1) (2014), 247--256.
\bibitem{michalak1} L. P. Michalak, \textsl{Realization of a graph as the Reeb graph of a Morse function on a manifold}. Topol. Methods in Nonlinear Anal. 52 (2) (2018), 749--762, arxiv:1805.06727.; doi: 10.12775/TMNA.2018.029. 
\bibitem{michalak2} L. P. Michalak, \textsl{Combinatorial modifications of Reeb graphs and the realization problem}, arxiv:1811.08031.; doi: 10.1007/s00454-020-00260-6.
\bibitem{milnor} J. Milnor, \textsl{Lectures on the h-cobordism theorem}, Math. Notes, Princeton Univ. Press, Princeton, N.J. 1965.
\bibitem{reeb} G. Reeb, \textsl{Sur les points singuliers d\'{}une forme de Pfaff compl\'{e}tement int\`{e}grable ou d\'{}une fonction num\'{e}rique}, Comptes Rendus
 Hebdomadaires des S\'{e}ances de I\'{}Acad\'{e}mie des Sciences 222 (1946), 847--849.
\bibitem{saeki} O. Saeki, \textsl{Reeb spaces of smooth functions on manifolds}, Intermational Mathematics Research Notices, maa301, https://doi.org/10.1093/imrn/maa301, arxiv:2006.01689.
\bibitem{sharko} V. Sharko, \textsl{About Kronrod-Reeb graph of a function on a manifold}, Methods of Functional Analysis and
 Topology 12 (2006), 389--396.
\end{thebibliography}
\end{document}